\documentclass[12pt]{amsart}
\usepackage{mathrsfs}
\usepackage{txfonts}
\usepackage{amssymb}

\makeatletter \@mparswitchfalse \makeatother
\def\eqref#1{(\ref{eq#1})}

\numberwithin{equation}{section}
\begin{document}
\title{Subdiagonal algebras with the Beurling type invariant subspaces}
\author{Guoxing Ji}
  \address{School  of Mathematics and Information Science,
  Shaanxi Normal University,
  Xian , 710062, People's  Republic of  China}
  \email{gxji@snnu.edu.cn}

     \thanks{This research was
supported by the National Natural
   Science Foundation of China(No. 11771261) and the Fundamental Research Funds for the Central Universities (Grant No. GK201801011).}
    %\no{\footnotesize{\bf MSC(2000):\hspace{2mm}Primary    46L52, 47L75: Secondary 46K50,  46J15 }

   % \subjclass{     46L52\sep  46K50}
 \maketitle
\begin{abstract}
Let $\mathfrak A$ be a maximal subdiagonal algebra in a $\sigma$-finite von Neumann algebra $\mathcal M$. If every right invariant subspace of $\mathfrak A$ in the non-commutative Hardy space $H^2$ is of Beurling type, then we say  $\mathfrak A$  to be type 1.  We determine   generators of these algebras and consider  a Riesz type factorization theorem for the  non-commutative $H^1$ space. We show that the right analytic Toeplitz algebra on the non-commutative Hardy space $H^p$ associated with  a type 1 subdiagonal algebra with multiplicity 1 is hereditary reflexive.\\
  {\bf keywords}{  von Neumann algebra, subdiagonal algebra,  non-commutative Hardy
space, factorization, reflexivity}
   \end{abstract}
\baselineskip18pt

\section{Introduction}

Beurling's  invariant subspace theorem plays a very important role in classical Hardy space theory and an  extensive version  of Beurling-Lax-Halmos theorem was  developed(cf.\cite{beu,hal,lax}). There are a lots of applications of this theorem since then.   Moreover, many commutative as well as noncommutative  extensions of this theorem have appeared over
the decades, for example, Weak$^*$- Dirichlet
algebras(\cite{sw}), non-self-adjoint crossed products(cf. \cite{mms} and references therein) and so on. On the other hand, Arveson in \cite{arv1} introduced the notion of subdiagonal algebras, as the noncommutative analogue of  the classical Hardy  space $H^{\infty}(\mathbb T)$, to unify several aspects of non-self-adjoint operator algebras. It is remarkable that there are several  successful noncommutative extensions of classical $H^p$ spaces based on  subdiagonal algebras(cf. \cite{ble,ble1,ble2,ble3,jio,jis,jig,jig1}). One important  extension  is due to  Blecher and    Labuschagnethe on  Beurling type invariant subspace  theorem for finite subdiagonal algebras(\cite{ble4}). Very recently, Labuschagne in \cite{lab} extend  their results to general maximal subdiagonal algebras in a $\sigma$-finite von Neumann algebra.  In fact, they decompose an invariant subspace as an internal $L^2$-column sum  of type 1 and type 2 invariant subspace  according to   $L^p$-column sums due to  Junge and Sherman \cite{js}. Moreover, every type 1 invariant subspace has the Beurling type, that is, every right invariant subspace $\mathfrak M$ in  the non-commutative $L^2$ space associated with  a $\sigma$-finite von Neumann algebra is an internal $L^2$-column sum of a family of the form $U_iH^2$, where $\{U_i:i\geq 1\}$ is a family of partial isometries such that $U_j^*U_i=0$ for $i\not=j$ while $U_i^*U_i(i\geq1)$ is a projection in the diagonal algebra of considered subdiagonal algebra in the  von Neumann algebra.  Note that every  invariant subspace in the classical Hardy space $H^2$ has  the Beurling  type.  It then becomes natural   when  every invariant subspace in the  non-commutative Hardy space $H^2$  has the Beurling  type. We consider those maximal subdiagonal algebras whose invariant subspaces have the  Beurling type in non-commutative $H^2$ space  in a $\sigma$-finite von Neumann algebra. We
   firstly
   recall some notions.

  Let $\mathcal M$ be
a $\sigma$-finite von Neumann algebra acting on a complex  Hilbert
$\mathcal H$. We denote by $\mathcal M_*$ the space of all
$\sigma$-weakly continuous linear functionals of $\mathcal M$. Let
$\Phi$  be  a faithful normal conditional expectation from $\mathcal
M$ onto  a von Neumann subalgebra   $\mathfrak D$. Arveson \cite{arv1} gave the following definition. A subalgebra
$\mathfrak A$ of $\mathcal M$, containing $\mathfrak D$, is called a
subdiagonal algebra of $\mathcal M$ with respect to $\Phi$ if

(i) $\mathfrak A \cap \mathfrak A^{*} =\mathfrak D$,

(ii) $\Phi$ is multiplicative on $\mathfrak A$, and

(iii)  $\mathfrak A + \mathfrak A^*$ is $\sigma$-weakly dense in
$\mathcal M$.

\noindent The algebra $\mathfrak D$ is called the diagonal of
$\mathfrak A$.
 Although subdiagonal  algebras are not assumed to be $\sigma$-weakly
 closed in  \cite{arv1},  the $\sigma $-weak closure of a subdiagonal algebra
 is  again  a subdiagonal algebra of  $\mathcal M$ with respect to
 $\Phi$(\cite[Remark 2.1.2]{arv1}). Thus we assume that our subdiagonal
 algebras are always $\sigma$-weakly closed.

  We say that
$\mathfrak A$ is a maximal subdiagonal algebra in $\mathcal M$ with
respect to $\Phi$ in case that $\mathfrak A$ is not properly
contained in any other subalgebra of $\mathcal M$ which is
subdiagonal with respect to $\Phi$. Put $\mathfrak A_0 =\{X \in
\mathfrak A: \Phi(X)=0 \}$ and  $\mathfrak A_m =\{X \in \mathcal M:
\Phi(AXB)=\Phi(BXA)=0,~ \forall A \in \mathfrak A,~ B \in \mathfrak
A_0 \}$.
 By   \cite[Theorem 2.2.1]{arv1}, we recall that $\mathfrak A_m $ is a maximal
 subdiagonal algebra of $\mathcal M$ with respect to $\Phi$ containing
 $\mathfrak A$.

 We  next recall Haagerup's noncommutative $L^p$ spaces associated with a
 $\sigma$-finite von Neumann algebra $\mathcal M$. Let $\varphi$ be a faithful
 normal state on $\mathcal M$ and  let $\{\sigma_t^{\varphi}:t\in\mathbb R\}$
  be the modular automorphism group of $\mathcal M$ associated with $\varphi$ by Tomita-Takesaki theory.
   We consider  the crossed product
 $\mathcal N=\mathcal M\rtimes_{\sigma^{\varphi}} \mathbb R$ of $\mathcal M$
 by $\mathbb R$ with respect to $\sigma^{\varphi}$.
Then we have that $\mathcal N $ is a von Neumann algebra on
$L^2(\mathbb R,\mathcal H)$ generated by the operators $\pi(x)$,
$x\in \mathcal M$, and $\lambda(s)$, $s \in \mathbb R$ defined by
the equations
  $$
(\pi(x)\xi)(t)=\sigma^{\varphi}_{-t}(x)\xi(t), \ \ \xi\in
L^2(\mathbb R,\mathcal H), \ t \in \mathbb R,
$$
and
$$
(\lambda(s)\xi)(t)=\xi(t-s),\ \ \xi\in L^2(\mathbb R,\mathcal H), \
t \in \mathbb R.$$
  We identify $\mathcal M$ with its image $\pi(\mathcal M)$ in
  $\mathcal N$.

  We denote by $\theta$ the dual action of $\mathbb R$ on $\mathcal
  N$. Then $\{\theta_s:s\in \mathbb R\}$  is an  automorphisms group
  of $\mathcal N$
  characterized by
  $\theta_s (X)=X, X\in\mathcal M$, $\theta_s (\lambda
  (t))=e^{ist},t\in\mathbb R$.

Note that $\mathcal M=\{X\in \mathcal N: \theta_s(X)=X,\forall
s\in\mathbb R\}$.  $\mathcal N$ is a semifinite von Neumann algebra
and there is the normal faithful semifinite trace $\tau$ on
$\mathcal N$ satisfying
$$\tau\circ \theta_s=e^{-s}\tau,  \ \ \  \forall s\in\mathbb R.$$

According to Haagerup \cite{haa,ter}, the noncommutative $L^p$
spaces $L^p(\mathcal M)$ for each $0< p\leq \infty$ is defined as
the set of all $\tau$-measurable  operators $x$   affiliated with
$\mathcal N$ satisfying
$$
\theta_s(x)=e^{-\frac{s}{p}}x,  \ \forall  s\in\mathbb R.$$ There is
a linear bijection  between the predual $\mathcal M_*$ of $\mathcal
M$ and $L^1(\mathcal M)$: $f\to h_f$. If we define
$\mbox{tr}(h_f)=f(I), f\in$$ \mathcal M_*$, then $$
\mbox{tr}(|h_f|)=\mbox{tr}(h_{|f|})=|f|(I)=\|f\|$$ for all
 $f\in \mathcal M_*$ and
$$|\mbox{tr}(x)|\leq\mbox{tr}(|x|)$$ for all  $x\in L^1(\mathcal M)$.
Note that for any $x\in L^p(\mathcal M)$,  $\|x\|_p=(tr(|x|^p))^{\frac1p}$ is the norm of $x$.
As in \cite{haa}, we define the operator $L_A$ and $R_A$ on
$L^p(\mathcal M)$($1\leq p<\infty)$  by $L_Ax=Ax$ and $R_Ax=xA$ for
all  $A\in\mathcal M$ and $x\in L^p(\mathcal M)$.
 Note that
$L^2(\mathcal M)$ is a Hilbert space with the inner product $\langle a,b \rangle
=\mbox{tr}(b^*a)$, $ \forall a,b \in L^2(\mathcal M)$ and
  $A\to L_A$( resp. $A\to R_A$) is a
  faithful representation (resp. anti-representation) of $\mathcal
  M$ on $L^2(\mathcal M)$. We may identify $\mathcal M$ with
  $L(\mathcal M)=\{L_A:A\in \mathcal M\}$ on $L^2(\mathcal M)$. We also denote by $R(\mathcal M)$ the right multiplication operators by $\mathcal M$, that is, $R(\mathcal M)=\{R_A:A\in\mathcal M\}$. It is well-known that  $L(\mathcal M)^{\prime}=R(\mathcal M)$ and $R(\mathcal M)^{\prime}=L(\mathcal M)$.

Let $h_0$ be the
noncommutative Radon-Nikodym  derivative of the dual weight of
$\varphi$ on $\mathcal N$ with respect to $\tau$.  Then $h_0$ is the
image($h_{\varphi}$) of $\varphi$ in $L^1(\mathcal M)$. For a subset $S\subseteq L^p(\mathcal M)$, denote by $\vee S=[S]_p$ the closed($\sigma$-weakly closed if $S\subseteq \mathcal M$) subspace  generated by $S$.  It is well-known that
$L^p(\mathcal M)=[h_0^{\frac{\theta}{p}}\mathcal Mh_0^{\frac{1-\theta}p}]_p$ for any $1\leq p<\infty$ and  $\theta\in[0,1]$.
Thus we may define
  the noncommutative $H^p$ space  $H^p(\mathcal M)$  and $H_0^p(\mathcal M)$ in $L^p(\mathcal M)$
for any $1\leq p<\infty$  as
$$H^p=H^p(\mathcal M)=[h_0^{\frac{\theta}{p}}\mathfrak Ah_0^{\frac{1-\theta}p}]_p\mbox{ and } H_0^p=H_0^p(\mathcal M)=[h_0^{\frac{\theta}{p}}\mathfrak A_0h_0^{\frac{1-\theta}p}]_p$$ for any $\theta\in[0,1]$(\cite[Definition
2.6]{jig}, \cite[Proposition 2.1]{jig1}).  It is known that the  noncommutative $H^p$   space associated
with a subdiagonal algebra $\mathfrak A$ with respect to the
faithful normal expectation $\Phi$  are  independent of the choice of
states which preserve   $\Phi$(\cite[Theorem 2.5]{jig}).
We also call  the maximal subdiagonal algebra $\mathfrak A$  as the  non-commutative $H^{\infty}$.

In this paper, we consider  maximal subdiagonal algebras, which we called type 1 subdiagonal algebras, whose right(resp. left) invariant subspaces in   the non-commutative Hardy space $H^2$ are of  type 1. That is, every right(resp. left) invariant subspace in $H^2$ has Beurling type.  We determine the generators of type 1 subdiagonal algebras in Section 1.  Moreover  we give a Riesz type factorization theorem for non-commutative $H^1$ which says that every element in $H^1$ is a product of two elements in $H^2$ in Section 2. In Section 3, we show that the right(resp. left) analytic Toeplitz algebra associated with  a type 1 subdiagonal algebra with multiplicity 1(defined in Section 2) on non-commutative $H^p$($1<p<\infty)$) is hereditary reflexive. This result generalized the Sarason's result on  classical analytic Toeplitz algebra(\cite{sara}) and Peligard's result on finite non-self-adjoint crossed products(\cite{pel}).

\section{Generators of a type 1 Subdiagonal algebra}

Let $\mathfrak A$ be a maximal subdiagonal algebra in  a von Neumann algebra $\mathcal M$ with respect to $\Phi$.
 We recall that a  closed subspace $\mathfrak M$
of $ L^2(\mathcal M)$  is right(resp. left)  invariant if  $\mathfrak M
\mathfrak A\subseteq \mathfrak M$(resp. $ \mathfrak A\mathfrak M
\subseteq \mathfrak M$).   By the symmetry, it is sufficient to consider one  side.  We next consider the right invariant subspaces. Following \cite{ble4,lab,nak}, we define the right wandering subspace of
$\mathfrak M$  to be the space $W = \mathfrak M\ominus [\mathfrak M\mathfrak A_0]_2$. We say that $\mathfrak M$ is type 1 if $W$ generates $\mathfrak M$ as an
$\mathfrak A$-module (that is, $\mathfrak M = [W\mathfrak A]_2$). We will say that $\mathfrak M$ is type 2 if $W = \{0\}$. Note that every   right invariant subspace $\mathfrak M$ is an $L^2$-column sum $\mathfrak M=\mathfrak N_1\oplus^{col}\mathfrak N_2$, where     $\mathfrak N_i$  is of  type$i$ for $i=1,2$ from \cite[Theorem 2.1]{ble4} and
\cite[Theorem 2.3]{lab}.  In particular,   if $\mathfrak M$ is    type 1, then $\mathfrak M$ has Beurling type, that is, there are a family  of partial isometries $\{U_n:n\geq 1\}$ satisfying $U_i^*U_j=0 $ for $i\not= j$ and $U_i^*U_i\in\mathfrak D$ such that $\mathfrak M=\oplus_{n}^{col}U_nH^2 $. We refer    to {\cite{ble4,lab} for  more details.

{\bf Definition 2.1.} Let $\mathfrak A$ be a  maximal subdiagonal algebra in $\mathcal M$. If every right invariant subspace of $\mathfrak A$ in $ H^2$ is  type 1, then we say that $\mathfrak A$ is a  type 1 subdiagonal algebra.

For any  positive integer $n\geq 1$, let  $\mathfrak A_0^n$ be the  $\sigma$-weakly closed ideal of $\mathfrak A$  generated by $\{a_1a_2\cdots a_n:a_j\in\mathfrak A_0\}$.   Put $\mathfrak M_n=[  H^2\mathfrak A_0^n]_2$. Then $\mathfrak M_n\subseteq H^2$ is  a right invariant subspace such that $\mathfrak M_{n+1}=[\mathfrak M_n\mathfrak A_0]_2$.

{\bf Proposition 2.1.} Let $\mathfrak A$ be a maximal subdiagonal algebra. Then  $\mathfrak A$ is of type 1 if and only if $\bigcap\limits_{n=1}^{\infty}\mathfrak M_n=\{0\}$.

%{\bf Proof}
 \begin{proof}$\Longrightarrow$  Note that
$[(\bigcap\limits_{n=1}^{\infty} \mathfrak M_n)\mathfrak A_0]_2=\bigcap\limits_{n=1}^{\infty} \mathfrak M_n$ is a type 2 right invariant subspace   in $H^2$. Then $\bigcap\limits_{n=1}^{\infty}\mathfrak M_n=\{0\}$.

$\Longleftarrow$ If $\mathfrak M\subseteq   H^2$ is a type 2 right invariant subspace, then
$[\mathfrak M \mathfrak A_0]_2=\mathfrak M$. Thus $\mathfrak M= [  H^2\mathfrak A_0^n]_2 =\mathfrak M_n$ for any $n$. It follows that  $\mathfrak M=\{0\}$ and  $\mathfrak A$ is of type 1.
\end{proof}
%\qed

  We recall that
  $\{\sigma_{t }^{\varphi}: t \in \mathbb R\}$ is the  modular automorphism
group of $\mathcal M$ associated with $\varphi$   and we have
that the following representation of
$\sigma_t^{\varphi}$(cf.\cite{kos1,kos2}):
\begin{equation}\sigma_t^{\varphi}(X)=h_0^{it}Xh_0^{-it},\ \forall t\in\mathbb
R, \ \forall X\in\mathcal M .
   \end{equation}
Since $\mathfrak A$ is $\{\sigma_t^{\varphi}:t\in\mathbb R\}$ invariant by \cite[Theorem 2.4]{jio},  $\mathfrak M_n=[ \mathfrak A_0^n H^2 ]_2$  for all $n$. This proposition  says that we may define type 1 subdiagonal algebras by use of left invariant subspaces. Thus type 1 subdiagonal algebras are independent of choice of right or left invariant subspaces. Moreover we also have the following property.
 We recall that $\mathfrak A$ has the universal factorization property if every invertible operator $T\in\mathcal M$ has a factorization $T=WA$, where $W\in\mathcal M$ is unitary and  $A,A^{-1}\in\mathfrak A$(\cite{pit}).

{\bf Corollary 2.2.} If $\mathfrak A$ is a type 1 subdiagonal allgebra, then  $\mathfrak A$ has the universal  factorization property.

 \begin{proof}
 Since  $\bigcap\limits_{n=1}^{\infty}\mathfrak M_n=\{0\}$, $\bigcap\limits_{n=1}^{\infty}\mathfrak A_0^n=\{0\}$. Then $\mathfrak A_0$ contains no nonzero idempotent. It follows from \cite[Theorem 3.8]{jis} that $\mathfrak A$ has the universal  factorization property. \end{proof}

Let  $\mathfrak A$ be  of  type 1. Note that $H_0^2\subseteq H^2$ is a right invariant subspace.  Then    there is a family of  partial isometries $\{U_n:n\geq 1\}$ in $\mathcal M$  such that   \begin{equation}
 H_0^2 = \oplus_{n\geq 1}^{col} U_n  H^2
 \end{equation}
 with $U_i^*U_j=0$  for $i\not=j$ and $U_j^*U_j\in \mathfrak D$. These partial isometries are not unique in general.  If    there are finite partial isometries such that formula $(2.2)$ holds, then  we define the multiplicity of $\mathfrak A$ as
 $$m(\mathfrak A)=\min\{|I|:  H_0^2=\oplus_{n\in I}^{col}U_n   H^2,  U_m^*U_n=0,  \mbox{ for } m\not=n \mbox{ and }U_n^*U_n\in \mathfrak D\},$$
 where $|I|$ is the number of elements of $I$. Otherwise we define the multiplicity of $\mathfrak A$ is $\infty$. We now fixed   such a family of partial isometries  $\mathcal U=\{U_n:n\geq 1\}$    in $\mathcal M$ such that the  formula $(2.2)$ holds.

 {\bf Proposition 2.3.} Let  $\mathfrak A$ be  a  type 1  subdiagonal algebra. Then
$\mathfrak A_0U_m^*\cup U_n^*\mathfrak A_0\subseteq \mathfrak A$,  $U_m\mathfrak A U_n^*\cup U_m^*\mathfrak A U_n\subseteq \mathfrak A$ and  $U_m\mathfrak DU_n^*\cup U_m^*\mathfrak DU_n\subseteq  \mathfrak D$ for all $m,n\geq 1$.

 \begin{proof}
 It is elementary  that $\mathfrak A_0=\{A\in\mathcal M: AH^2\subseteq H_0^2\} $ from \cite[Theorem 2.7]{jig}. Note that    $U_nH^2\subseteq H_0^2$ for all $n\geq 1$. Then $U_n\in \mathfrak A_0$.

  Now    $U_n^* H_0^2=U_n^*U_nH^2\subseteq H^2$ for any $n\geq1$. Let $A\in\mathfrak A_0$. Then $A U_n^*H_0^2=AU_n^*U_nH^2\subseteq AH^2\subseteq H_0^2$. It follows  from \cite[Theorem 2.7 ]{jig}  that $A U_n^*\in \mathfrak A$. Furthermore,   $U_m^*AH^2\subseteq U_m^*H_0^2=U_m^*U_mH^2 \subseteq H^2$.  This means again that $U_m^*A\in\mathfrak A$.

   On the other hand,  let $B\in\mathfrak A$.  $B U_n^*H_0^2=BU_n^*U_nH^2\subseteq BH^2\subseteq H^2$ for all  $n\geq 1$. Then $U_mB U_n^*H_0^2\subseteq H_0^2$ for all     $m,n\geq 1$.
  It follows  that
 $U_mBU_n^*\in\mathfrak A$ for all   $m,n\geq 1$   by \cite[Theorem 2.7 ]{jig}. Again,
 $U_m^*BU_nH^2\subseteq U_m^*BH_0^2\subseteq U_m^*H_0^2=U_m^*U_mH^2\subseteq H^2$ for any $m,n\geq 1$.  Thus $U_m^*BU_n\in\mathfrak A$.

   Take any $D\in\mathfrak D$. $U_mDU_n^*,U_m^*DU_n\in\mathfrak A$. $(U_mDU_n^*)^*=U_nD^*U_m^*\in\mathfrak A$ and $(U_m^*DU_n)^*=U_n^*D^*U_m\in\mathfrak A$.  Then $U_mDU_n^*, U_m^*DU_n\in\mathfrak D$ for all $m,n\geq 1$.\end{proof}

{\bf Corollary 2.3. }  Let $E_n=U_n^*U_n$, $F_n=U_nU_n^*$ for all $n\geq 1$.  Then  $U_n^*\mathfrak A U_m=E_n\mathfrak AE_m$ and $U_n\mathfrak AU_m^*=F_n\mathfrak AF_m$ for all $m,n\geq 1$.

 \begin{proof} $E_n\mathfrak AE_m=U_n^*U_n\mathfrak AU_m^*U_m\subseteq U_n^*\mathfrak A U_m=E_nU_n^*\mathfrak A U_mE_m\subseteq E_n\mathfrak AE_m$. Similarly we have that $U_n\mathfrak AU_m^*=F_n\mathfrak AF_m$.\end{proof}

 Put $W_0=L^2(\mathfrak D)$ and $W_n=\mathfrak M_n\ominus [\mathfrak M_n\mathfrak A_0]_2=\mathfrak M_n\ominus \mathfrak M_{n+1}$ is the wandering subspace for $\mathfrak M_n$ for $n\geq 1$. Then $H^2=\bigoplus\limits_{n=0}^{\infty}W_n$ by Proposition 2.1.

{\bf Proposition 2.4. } Let $\mathfrak A$ be a type 1 subdiagonal algebra. Then for any $n\geq 1$,  $\mathfrak M_n =\vee\{U_{i_1}U_{i_2}\cdots U_{i_n}H^2: U_{i_k}\in\mathcal U\} $ and
$W_n=\vee\{U_{i_1}U_{i_2}\cdots U_{i_n}L^2(\mathfrak D): U_{i_k}\in\mathcal U \}$.

 \begin{proof}
 $H^2=L^2(\mathfrak D)\oplus H_0^2$. Then  $\mathfrak M_1=H_0^2=\oplus_{n\geq 1}^{col}U_nH^2=\vee\{U_nL^2(\mathfrak D):n\geq 1\}\oplus \vee\{U_nH_0^2:n\geq 1\}$. This  means that $W_1=\vee\{U_nL^2(\mathfrak D):n\geq 1\}$. Thus the conclusion holds for $n=1$.
             Assume that the conclusion holds for all  $k\leq n-1$.
               Then
 \begin{align*}
 \mathfrak M_n &=[\mathfrak M_{n-1}\mathfrak A_0]_2\\
 &=[\vee\{U_{i_1}U_{i_2}\cdots U_{i_{n-1}}H^2:U_{i_k}\in\mathcal U\} \mathfrak A_0]_2\\
 &=\vee\{U_{i_1}U_{i_2}\cdots U_{i_{n-1}}[H^2\mathfrak A_0]_2:U_{i_k}\in\mathcal U\} \\
 &=\vee\{U_{i_1}U_{i_2}\cdots U_{i_{n-1}}H^2_0:U_{i_k}\in\mathcal U\} \\
 &=\vee\{U_{i_1}U_{i_2}\cdots U_{i_{n-1}}(\oplus_{n\geq1}^{col}U_nH^2):U_{i_k}\in\mathcal U\} \\
 &=\vee\{U_{i_1}U_{i_2}\cdots U_{i_{n}}H^2:U_{i_k}\in\mathcal U\} .
  \end{align*}
  On the other hand, $W_n=\mathfrak M_n\ominus [\mathfrak M_n\mathfrak A_0]=\mathfrak M_n \ominus \mathfrak M_{n+1}$. It is known that
  $$(U_{j_1}U_{j_2}\cdots U_{j_n})^*(U_{i_1}U_{i_2}\cdots U_{i_n})\in\mathfrak D$$ for all $i_k,j_k\geq1$ from Proposition 2.3.
   It easily follows that $\vee\{U_{i_1}U_{i_2}\cdots U_{i_n}L^2(\mathfrak D): U_{i_k}\in\mathcal U\}\subseteq W_n$.
   However,
   \begin{align*}
    \mathfrak M_n &=\vee\{U_{i_1}U_{i_2}\cdots U_{i_n}H^2: U_{i_k}\in\mathcal U\}\\
    &=\vee\{U_{i_1}U_{i_2}\cdots U_{i_n}(L^2(\mathfrak D)\oplus H_0^2): U_{i_k}\in\mathcal U\}\\
    &= \vee\{U_{i_1}U_{i_2}\cdots U_{i_n}L^2(\mathfrak D): U_{i_k}\in\mathcal U\}\oplus \mathfrak M_{n+1}.
    \end{align*}
    Thus
 $W_n=\vee\{U_{i_1}U_{i_2}\cdots U_{i_n}L^2(\mathfrak D): U_{i_k}\in\mathcal U\}$.\end{proof}

We now give generators of a type 1 subdiagonal algebra.

{\bf Theorem 2.5.} Let $\mathfrak A$ be a  type 1 subdiagonal algebra. Then  $\mathfrak A_0=\vee\{U_{i_1}U_{i_2}\cdots U_{i_n}\mathfrak D:U_{i_k}\in\mathcal U,  n\geq 1\}$ and $\mathfrak A=\mathfrak D+\mathfrak A_0$.

   \begin{proof}
   Put $\mathcal A_0=\vee\{U_{i_1}U_{i_2}\cdots U_{i_n}\mathfrak D:U_{i_k}\in\mathcal U, \ n\geq 1\}$. It is trivial that $\mathcal A_0\subseteq\mathfrak A_0$. We note that $D\mathfrak M_n=\mathfrak M_n$ for all $n\geq 1$. Then $\mathfrak DW_n=W_n$ for all $n\geq1$. Thus $\mathfrak D(\oplus_{n\geq 1}U_nL^2(\mathfrak D))=\oplus_{n\geq 1}U_nL^2(\mathfrak D)$. Take any $D\in\mathfrak D$ and $m\geq 1$.
Then $DU_mh_0^{\frac12}=\oplus_{n\geq 1}U_n\xi_n$ for some $\xi_n\in L^2(\mathfrak D)$.  Thus $U_nU_n^*DU_mh_0^{\frac12}=
U_nU^*_nU_n\xi_n=U_n\xi_n$ for all $n\geq 1$. That is,
\begin{equation}DU_mh_0^{\frac12}=\oplus_{n\geq 1}U_nU_n^*DU_mh_0^{\frac12}.\end{equation}
We know that $U_n^*DU_m=D_{nm}\in\mathfrak D$  for all $n\geq 1$ by Proposition 2.3. Put $A_{nm}=\sum\limits_{k=1}^nU_kD_{km}$. Then \begin{align*}
A_{nm}^*A_{nm}&=(\sum\limits_{k=1}^nU_kD_{km})^*(\sum\limits_{k=1}^nU_kD_{km})
=\sum\limits_{k=1}^nD_{km}^*U_k^*U_kD_{km}\\
&=\sum\limits_{k=1}^nU_m^*D^*U_kU_k^*DU_m
=U_m^*D^*(\sum\limits_{k=1}^nU_kU_k^*)DU_m\\
&\leq U_m^*D^*DU_m.
\end{align*}
   That is, $\{A_{nm}:n\geq 1\}$ is a bounded sequence in $\mathcal M$. It then follows from $(2.3)$ that $A_{nm}$ converges $\sigma$-weakly and
   \begin{equation}DU_m=\lim\limits_{n\to\infty}A_{nm}=\sum\limits_{n=1}^{\infty}U_nD_{nm}.\end{equation}
    Thus $DU_m\in\mathcal A_0$ for all $D\in\mathfrak D$ and $m\geq1$. This implies that  $\mathcal A_0$ is a subalgebra of $\mathfrak A_0$.

     We claim that $\mathcal A_0+\mathfrak D+\mathcal A_0^*$ is a $^*$-subalgebra of $\mathcal M$.  It is elementary that $\mathcal A_0$ is a $\mathfrak D$ bimodule.  For any  $ D_1, D_2\in \mathfrak D$ and $n,m\geq 1$, put
     $A=D_1(U_{i_1}U_{i_2}\cdots U_{i_n})^*$ and $B=(U_{j_1}U_{j_2}\cdots U_{j_m})D_2$. Then $AB, BA\in \mathcal A_0$ if $n<m$, $AB, BA\in\mathfrak D$ if $n=m$ and $AB, BA\in \mathcal A_0^*$ if $n>m$ by Proposition 2.3. Thus  $\mathcal A_0+\mathfrak D+\mathcal A_0^*$ is a $^*$-subalgebra of $\mathcal M$.

     It is easily shown  that  $(\mathcal A_0+\mathfrak D+\mathcal A_0^*)h_0^{\frac12}$ is dense in $ L^2(\mathcal M)$.
     Let  $\mathcal M_1$ be the  $\sigma$-weak closure of $\mathcal A_0+\mathfrak D+\mathcal A_0^*$. Then $\mathcal M_1$ is a von Neumann subalgebra of $\mathcal M$ and $\mathcal A=\mathcal A_0+\mathfrak D$ is  a subdiagonal algebra  of $\mathcal M_1$ with respect to $\Phi$ such that $L^2(\mathcal M_1)= L^2(\mathcal M)$.   It is easily known  that $\mathcal M_1=\mathcal M$ since $h_0^{\frac12} $ is a common cyclic and separating vector for $\mathcal M_1$ and $\mathcal M$. Thus $\mathcal A$ is a   subdiagonal algebra of $\mathcal M$ with respect to $\Phi$ such that  $\mathcal A\subseteq \mathfrak A$.
 If we may prove that $\mathcal A_0$ is $\{\sigma_t^{\varphi}:t\in\mathbb R\}$ invariant, then we have  that $\mathcal A$ is maximal subdiagonal by  \cite[Theorem 1.1]{xu} and thus $\mathcal A_0=\mathfrak A_0$.

   Since $\mathfrak A_0$  and $\mathfrak D$ are   $\{\sigma_t^{\varphi}:t\in\mathbb R\}$ invariant from \cite[Theorem 2.4]{jio},
    \begin{equation}H_0^2=\oplus_{n\geq 1}^{col}U_nH^2=\oplus_{n\geq 1}^{col} \sigma_t^{\varphi}(U_n)H^2
     \end{equation}
     for all $t\in\mathbb R$. It follows  from above formula that
    $$U_m^*\sigma_t^{\varphi}(U_n)H^2\subseteq U_m^*H_0^2= U_m^*U_mH^2\subseteq H^2$$
    and $$
    (\sigma_t^{\varphi}(U_n))^*U_mH^2\subseteq (\sigma_t^{\varphi}(U_n))^*H_0^2= (\sigma_t^{\varphi}(U_n))^*\sigma_t^{\varphi}(U_n)H^2=\sigma_t^{\varphi}(U_n^*U_n)H^2\subseteq H^2,$$
     which imply that
     \begin{equation}
     D_{nm}(t)=U_m^*\sigma_t^{\varphi}(U_n)\in\mathfrak D\end{equation}
      for all $m,n\geq 1$ and $t\in\mathbb R$  by \cite[Theorem 2.7 ]{jig}.

     Now for any $m\geq 1$  and $t\in\mathbb R$,
     $\sigma_t^{\varphi} (U_m)h_0^{\frac12}\in H_0^2=\oplus_{n\geq1}^{col}U_nH^2$. Then
     $\sigma_t^{\varphi} (U_m)h_0^{\frac12}=\oplus_{n\geq 1}U_n\xi_n$ for some $\xi_n\in H^2$. In fact,
      $U_n\xi_n=U_nU_n^*U_n\xi_n=U_nU_n^*\sigma_t^{\varphi} (U_m)h_0^{\frac12}$.
   Thus we similarly have
   \begin{equation}\sigma_t^{\varphi}(U_m)=\sum_{n\geq 1}U_n(U_n^*\sigma_t^{\varphi}(U_m))=\sum_{n\geq 1}U_nD_{nm}(t)\end{equation} $\sigma$-weakly.
It follows that $\sigma_t^{\varphi}(U_m)\in\mathcal A_0$. Therefore $\mathcal A_0 $ is $\{\sigma_t^{\varphi}:t\in\mathbb R\}$ invariant and $\mathcal A_0=\mathfrak A_0$. Consequently, $\mathfrak A=\mathfrak D+\mathcal A_0$. \end{proof}

We recall that a subdiagonal algebra is anti-symmetric if $\mathfrak D=\mathbb C I$(\cite{arv1}).

{\bf Corollary 2.6.} Let $\mathfrak A $ be an anti-symmetric  subdiagonal algebra of type 1. Then $\mathcal M$ is $*$-isomorphic to $L^{\infty}(\mathbb T)$ and $\mathfrak A$ is isometrically isomorphic  as well as $\sigma$-weakly homeomorphic to $H^{\infty}(\mathbb T)$.

\begin{proof}
If $\mathfrak D=\mathbb C I$, then  $L^2(\mathfrak D)=\mathbb C$ and $H_0^2= UH^2$ for a unitary operator  $U\in\mathfrak A_0$ since $U^*U, UU^*\in \mathfrak D$ by Proposition 2.3. It follows that $\mathfrak A=\bigvee\{U^n:n\geq 0\}$ and $\mathcal M=\bigvee\{U^n: n\in\mathbb Z\}$ by Theorem 2.6. This means that $L_U $ is a bilateral shift acting on $L^2(\mathcal M)$ by Proposition 2.4. Thus
   $$H^2= \oplus_{n\geq 0}U^nL^2(\mathfrak D)=\oplus_{n\geq 0}U^n\mathbb C.$$
   Let $\rho(U)=M_z$, where $M_z $ is the multiplication operator by $z$ on $L^2(\mathbb T)$. Then $\rho$ may be  extended  a $*$-isomorphism from
   $\mathcal M$ onto $L^{\infty}(\mathbb T)$. It is trivial that $\rho(\mathfrak A)=H^{\infty}(\mathbb T)$ and $\rho$ is also a $\sigma$-weak homeomorphism.\end{proof}

{\bf Corollary 2.7.} Let $\mathfrak A$  be a type 1 subdiagonal algebra  and   $\mathfrak M\subseteq  L^2(\mathcal M)$
 a type 2 right invariant subspace. Then $\mathfrak M$ is right reducing and $\mathfrak M=EL^2(\mathcal M)$ for some  projection $E\in\mathcal M$.

  \begin{proof}
  It is known that $\mathfrak M$ is  $\mathfrak D$ right  reducing.  Note that  $[\mathfrak M\mathfrak A_0]_2=\mathfrak M$.  Since $\vee\{\mathfrak MU_n:n\geq1\}\supseteq \vee\{\mathfrak MU_{i_1}U_{i_2}\cdots U_{i_n}:i_k\geq 1,n\geq1\}$, $\mathfrak M=\vee\{\mathfrak M U_n:n\geq1\}$ by Theorem 2.5. For any $m,n \geq 1$ and $x\in \mathfrak M$, $R_{U_m^*}R_{U_n}x=xU_nU_m^*\in \mathfrak M$ since $U_nU_m^*\in\mathfrak D$ from Proposition 2.3.   Then $\mathfrak M$ is reducing.  Thus the projection onto $\mathfrak M$ is in the commutant $L(\mathcal M)$ of  $R(\mathcal M)$. Then  $\mathfrak M=EL^2(\mathcal M)$ for some  projection $E\in\mathcal M$.\end{proof}

 We now give two examples of type 1 subdiagonal algebras with multiplicity 1.
 Let $\mathcal H$ be a separable complex Hilbert space. $\mathcal E=\{e_n: 1\leq n\leq d\}$ is an orthonormal  basis, where $d\leq \infty$.
Put $\mathfrak D$ be the algebra of all bounded diagonal operators with respect to $\mathcal E$. Denote by $E_n$ the projection onto $[e_n]$ for all $n\geq 1$ and $\Phi(T)=\sum\limits_{n=1}^{d}E_nTE_n$, $\forall T\in\mathcal B(\mathcal H)$. Then $\Phi$ is a faithful normal conditional expectation from $\mathcal B(\mathcal H)$ onto $\mathfrak D$. Let $U$ be the unilateral  shift operator, $Ue_n=e_{n+1}$, $\forall n\geq 1$ if $\dim \mathcal H=\infty$ and $Ue_n=e_{n+1}$ for $n<d$ and  $Ue_d=0$ if $\dim\mathcal H=d<\infty$.

  {\bf Example 2.8.} Let $\mathfrak A=\bigvee\{U^n\mathfrak D: n\geq 0\}$. Then $\mathfrak A$ is a type 1 subdiagonal algebra of $\mathcal B(\mathcal H)$  with multiplicity 1.

  Another example is  the non-self-adjoint crossed product defined in \cite{mms}. Let $\mathfrak D$ be a von Neumann algebra on $\mathcal H$  and $\alpha$ an automorphism of $\mathfrak D$.
  We  consider  the crossed product
 $\mathfrak D\rtimes_{\alpha} \mathbb Z$ of $\mathfrak D$  by $\alpha $.
Then we have that $\mathfrak D\rtimes_{\alpha}\mathbb Z $ is the  von
Neumann algebra on $\ell^2(\mathbb Z,\mathcal H)$ generated by the operators
$\pi(D)$, $D\in \mathfrak D$, and $U$ defined by
the equations
  $$
(\pi(D)\xi)_n=\alpha^{-n}(D)\xi_n, \ \ \forall  \xi=\{\xi_n\}\in
\ell^2(\mathbb Z,\mathcal H),
$$
and
$$
(U\xi)_n=\xi_{n-1},\ \ \forall \xi=\{\xi_n\}\in \ell^2(\mathbb Z,\mathcal H).$$

We identify $\pi(\mathfrak D)$ with $\mathfrak D$. We have the following example by \cite[Proposition 3.5]{mms}.

{\bf Example 2.9.} Let $\mathfrak D\rtimes_{\alpha} \mathbb Z_+=\bigvee\{U^n\pi(\mathfrak D):n\geq 0\}$  be the non-self-adjoint  crossed product.  Then $\mathfrak D\rtimes_{\alpha} \mathbb Z_+$ is a type 1 subdiagonal algebra of $\mathfrak D\rtimes_{\alpha} \mathbb Z$  with multiplicity 1.

\section{Riesz type factorization for type 1 subdiagonal algebras}

 In the  classical function space $H^1(\mathbb T)$, we may factor a function $h\in H^1(\mathbb T)$ as a product of two functions $ h_ 1$ and $h_2$  in $H^2(\mathbb T)$.  Marsalli and  West in \cite{mar}  gave an analogue factorization theorem for  the non-commutative  Hardy space $H^1$ associated with  a finite subdiagonal algebra. However,  whether the factorization theorem holds  for a maximal subdiagonal algebra in  a $\sigma$-finite von Neumann algebra is unknown.  Very recently, Bekjan in \cite{bek} considered some related problems. We  consider a Riesz type factorization for type 1 subdiagonal algebras in this section.
 The following lemmas  are elementary. For completeness, we give  simple proofs.

 Next lemma holds for any maximal subdiagonal algebra.  We denote by $\mathcal M_p$ the set of all projections in $\mathcal M$.

{\bf Lemma 3.1. } Let $\mathfrak A$ be a   maximal subdiagonal algebra. Put $Q=\vee\{P\in\mathcal M_p:    P\mathcal M\subseteq \mathfrak A\}$. Then $(I-Q)\mathfrak AQ=\{0\}$.

\begin{proof}  By  assumption, $Q\mathcal M\subseteq \mathfrak A$. Then $Q\mathcal MQ\subseteq \mathfrak D$.
Thus $Q\mathfrak A_0^*\subseteq \mathfrak A \cap \mathfrak A_0^*=\{0$\} since $Q\in\mathfrak D$. It follows that
$(I-Q)\mathfrak AQ=(I-Q)\mathfrak DQ$. Note that $(I-Q)\mathfrak DQ\mathcal M\subseteq \mathfrak A$. Put
$\mathfrak M=[ (I-Q)\mathfrak DQ\mathcal Mh_0^{\frac12}]_2$. Then $\mathfrak M\subseteq H^2$ is a right reducing subspace since $\mathfrak A$ is $\{\sigma_t^{\varphi}:t\in\mathbb R\}$ invariant by \cite[Theorem 2.4]{jio}.  Then there is a projection $E\in\mathcal M $ such that $\mathfrak M=E L^2(\mathcal M)\subseteq H^2$.
This means that $E\mathcal M\subseteq \mathfrak A$. It is trivial that $E\leq I-Q$.  Then $E=0$  and $(I-Q)\mathfrak AQ=\{0\}$.\end{proof}

 {\bf Lemma 3.2. } Let $\xi\in L^2(\mathcal M)  $ be a right wandering vector of $\mathfrak A$ and $\xi=U |\xi|$ the polar decomposition of $\xi$. Then $U^*U\in\mathfrak D$ and $|\xi|\in L^2(\mathfrak D).$

 \begin{proof}
  It is known  that $|\xi|^2\in L^1(\mathfrak D)$(cf.\cite[Remark 2.1]{ble4}, \cite[Theorem 2.3]{lab}). Thus
 $|\xi|\in L^2(\mathfrak D)$ by \cite[Lemma 2.4]{jig}. On the other hand, $L_{U^*U}$ is the projection onto $[|\xi|\mathcal M]_2$. It is trivial that $[|\xi|\mathcal M]_2=[|\xi|\mathfrak A_0^*]_2\oplus [|\xi|\mathfrak D]_2\oplus [|\xi|\mathfrak A_0]_2$. This means that $U^*U|\xi|D=|\xi|D\in L^2(\mathfrak D)$ for all $D\in\mathfrak D$. Now for any $\eta\in L^2(\mathfrak D)\ominus [|\xi|\mathfrak D]_2$, $\langle \eta, x+y^*\rangle =tr((x^*+y)\eta)=0$ for all $x,y\in H_0^2$. Note that $[|\xi|\mathfrak A_0^*]_2\subseteq (H_0^2)^*$ as well as
$[|\xi|\mathfrak A_0]_2\subseteq H_0^2$. It follows that $U^*U\eta=0$. Thus $U^*U L^2(\mathfrak D)\subseteq L^2(\mathfrak D)$ and $U^*U\in\mathfrak D$. \end{proof}

{\bf Lemma 3.3.  } Let  $\mathfrak A$ be  a  type 1 subdiagonal algebra and $x\in L^2(\mathcal M)$ a nonzero vector.  Put $\mathfrak M=[x\mathfrak A]_2$.  If $\mathfrak M=\mathfrak N_1\oplus^{col}\mathfrak N_2$ such that $\mathfrak N_i$ is of type {i} for $i=1,\ 2$, then there is a partial isometry $U\in\mathcal M$ with  $U^*U\in\mathfrak D$ such that $\mathfrak N_1=UH^2$.

  \begin{proof} Put $W=[x\mathfrak A]_2\ominus [x\mathfrak A_0]_2$ and $x=x_1+x_2$ for $x_1\in W$ as well as $x_2\in[x\mathfrak A_0]$.
 It is sufficient to show that $W=[x_1\mathfrak D]_2$. Note that  $[x_1\mathfrak A_0]_2,[x_2\mathfrak D]_2\subseteq [x\mathfrak A_0]_2$ and $[x_1\mathfrak A]_2=[x_1\mathfrak D]_2\oplus [x_1\mathfrak A_0]_2$. Then $[x\mathfrak A]_2=[x_1\mathfrak D]_2\oplus [x\mathfrak A_0]_2$.
  That is, $W=[x_1\mathfrak D]_2$. The desired result follows from \cite[Proposition 2.4]{lab} and Lemma 2.2. \end{proof}

   Note that a vector $x\in L^2(\mathcal M)$ is right separating(resp. cyclic) if it is a separating(resp. cyclic) vector of $R(\mathcal M)$. We recall that  an element $\xi\in H^2$ is right outer if $[\xi \mathfrak A]=H^2$(cf.\cite{ble2}).

{\bf Lemma 3.4.} Let $x,y \in L^2(\mathcal M)$ be   right separating and cyclic vectors  such that $yx=h_0 $. Then there is a unitary operator  $U\in\mathcal M$ such that $x=U\xi $ for some right outer element  $\xi\in H^2$.

  \begin{proof} Put $\mathfrak M=[x\mathfrak A]_2$. Then $\mathfrak M$ is a right invariant subspace.  It is not right reducing. In fact if it is, then $[x\mathfrak A]_2=[x\mathcal M]_2$ and therefore $H^1=[h_0 \mathfrak A]_1=[yx\mathfrak A]_1=[yx\mathcal M]_1= L^1(\mathcal M)$, a contradiction. It now follows that $[x\mathfrak A]_2=UH^2 \oplus ^{col} E L^2(\mathcal M)$ for a partial isometry $U$ and a projection $E\in\mathcal M$ such that $E L^2(\mathcal M)$ is of type 2 by \cite[Theorem2.3]{lab} and Lemma 3.3. Thus $[yE L^2(\mathcal M)]_1\subseteq [yx\mathfrak A]_1=H^1$ is a right  invariant subspace in $H^1$ such that  $ [yE L^2(\mathcal M)\mathfrak A_0]_1=[yE L^2(\mathcal M)]_1$ since $[E L^2(\mathcal M)\mathfrak A_0]_2=E L^2(\mathcal M)$.
   It is known that $[yE L^2(\mathcal M)]_1$ is a right $\mathcal M$ invariant subspace. By \cite[Theorem 2.7]{tak},  there is a projection $F\in\mathcal M$ such that $[yE L^2(\mathcal M)]_1=F L^1(\mathcal M).$ Thus $F L^1(\mathcal M)\subseteq H^1$, which implies that $F\mathcal M\subseteq \mathfrak A$ by \cite[Theorem 2.7]{jig}. If $F\not=0$, then by Lemma 3.1,  we have $F L^1(\mathcal M)F\subseteq L^1(\mathfrak D)$. This implies that    $  [yE L^2(\mathcal M)\mathfrak A_0]_1\not=[yE L^2(\mathcal M)]_1$, a contradiction. Hence $F=0$. Thus $yE=0$. This implies that $E=0$ and $\mathfrak M=UH^2$.

Let $\xi \in H^2$ such that $x=U\xi$. Then $[x\mathcal M]_2\subseteq  U L^2(\mathcal M)$, which implies that $U$ is surjective.
On the other hand,
$x$ is a separating and cyclic vector for $\mathcal M$. We have that $x_1$  defined in Lemma 3.3 is right separating. In fact,
 for any $D\in \mathfrak D$, if  $x_1D=0$, then $xD=x_2D$ and  $h_0D=yxD=yx_2D\in [y[x\mathfrak A_0]_2]_1=[yx\mathfrak A_0]_1=H^1_0$. Therefore  $D=0$.  Now take any $B\in\mathcal M$ such that $x_1B=0$. Note that $x_1 $ is a wandering vector. Then $\langle x_1B,x_1B \rangle=tr(|x_1|^2BB^*)=tr(|x_1|^2\Phi(BB^*))=\langle x_1D, x_1D\rangle=0$, where $D=(\Phi(BB^*))^{\frac12}\in\mathfrak D$. Then $D=0$ and $B=0$. That is , $x_1 $ is right separating. It follows  that  $|x_1|$ is both   right separating and right cyclic in this case. Note that  $U|x_1|B=x_1B$ for all $B\in\mathcal M$. Then $U$ is injective. Thus $U$ is unitary.
 $UH^2= [x\mathfrak A]_2=[U\xi\mathfrak A]_2=U[\xi\mathfrak A]_2$. Then $H^2=[\xi\mathfrak A]_2$ and $\xi$ is a right outer  element.\end{proof}

 We now give main  results of this section.

{\bf Theorem 3.5.  } Let $\mathfrak A$ be a type 1 subdiagonal algebra and  $\eta \in L^2(\mathcal M)$  a nonzero vector. Then for any $\varepsilon>0$, there are a contraction $B\in\mathcal M$ and a right outer element  $\xi\in H^2$  with $\|\xi\|_2<\|\eta \|_2+\varepsilon $ such that $\eta =B\xi$. If $\eta\in H^2$, then we may have $B\in\mathfrak A$.

 \begin{proof} Put $x=(\eta^*\eta+\varepsilon^2h_0)^{\frac12}$. Then we have that $x\in  L^2(\mathcal M)$ such that $\varepsilon^2h_0\leq x^2 $, $\eta^*\eta\leq x^2 $ and  $ \|x\|_2^2=\|\eta \|_2^2+\varepsilon^2$.  By \cite[Formula (1.1)]{js}, there are injective contractions $A, C\in\mathcal M$ with dense ranges such that $\varepsilon^{-1}Ax=h_0^{\frac12}$ and $\eta=Cx$. By Lemma 3.4, there are a unitary operator $U\in\mathcal M$ and  a right outer element  $\xi\in H^2$
 such that $x=U\xi$. Thus $\eta=CU\xi$. Put $B=CU$. We have  $\eta=B\xi$ and $\|\xi\|_2=\|x\|_2<\|\eta\|_2+\varepsilon$.

 If  $\eta \in H^2$, then $B\xi A=\eta A\in H^2$ for all $A\in \mathfrak A$. Thus $B H^2\subseteq H^2$ and  $B\in \mathfrak A$.
  \end{proof}

 {\bf Theorem 3.6. } Let $\mathfrak A$ be a type 1 subdiagonal algebra and     $h\in H^1$. Then for any $\varepsilon>0$,
  there are $h_1,h_2\in H^2$  with $\|h_i\|_2<\|h\|_1^{\frac12}+\varepsilon$  for $i=1,2$   such that $h=h_1h_2$.
  If $h\in H_0^1$, we may choose  one of $h_i$ in $H_0^2$.

   \begin{proof}
   Let $h=V|h|=V|h|^{\frac{1}2}|h|^{\frac{1}2}$. Then $|h|^{\frac{1}2}\in  L^2(\mathcal M)$ with $\||h|^{\frac{1}2}\|_2=\|h\|_1^{\frac{1}2}$.  By Theorem 3.5, there are a contraction $B\in \mathcal M$ and a right  outer element $\xi\in H^2$  with $\|\xi\|_2<\||h|^{\frac12}\|_2+\varepsilon=\|h\|_1^{\frac12}+\varepsilon$ such that $|h|^{\frac{1}2}=B\xi$. Put $h_1=V|h|^{\frac12}B$, and $h_2=\xi$. We show that $h_1\in H^2$.
   In fact, $tr(hA)=0$ for all $A\in \mathfrak A_0$ since $h\in H^1$. Then
   $tr(h A)=tr(h_1h_2A)=0$ for all  $A\in \mathfrak A_0$.  Note that $h_2$ is outer. Then $[h_2\mathfrak A_0]_2=H^2_0$.
    Therefor  $h_1\in H^2$.  $\|h_1\|_2\leq \||h|^{\frac12}\|_2< \|h\|_1^{\frac12}+\varepsilon$. Moreover, if $h\in H_0^1$, then $tr(hA)=0$ for all $A\in \mathfrak A$. It easily follows that
   $h_1\in H_0^2$  in this case.\end{proof}

    Let $1\leq p,q\leq \infty$  such that $\frac1p +\frac1q=1$ and  $h\in H^1$. Does $h=h_1h_2$ for some $h_1\in H^p$ and $h_2\in H^q$?

\section{Analytic Toeplitz algebras for type 1 subdiagonal algebras }

In this section, we consider  the reflexivity of analytic Toeplitz algebras associated with a type 1 subdiagonal algebra. Let $X$ be a Banach space and $\mathcal B(X)$ the algebra of all bounded linear operators on $X$. A subspace $S\subseteq \mathcal B(X)$ is said to be reflexive if $S=\{A\in \mathcal B(X): Ax\in [Sx]\mbox{ for every }x\in X\}$. $S$ is called hereditary reflexive if $S$ is reflexive and every weak operator topology closed subspace of $S$ is reflexive(\cite{log}).  If  $S$ is an algebra containing the identity $I$ and denote by $Lat S$ the invariant subspace lattice of $S$ in $X$, then $S$ is reflexive if $S=algLat S=\{A\in \mathcal B(X): AM\subseteq M,\forall M\in Lat S\}$. Sarason \cite{sara} proved that the algebra of all analytic Toeplitz operators on classical Hardy space $H^2(\mathbb T)$ is hereditary reflexive. Peligrad \cite{pel} extend to the case of non-commutative Hardy spaces $H^p$ for $1<p<\infty$  associated with  finite crossed products. We may extend to those for  type 1 subdiagonal algebras. In fact, the key step is to determine  when the  left and right analytic Toeplitz algebra associated with a  maximal subdiagonal algebra on the non-commutative Hardy space $H^p$   are commutants of each other. We have done when  $p=2$ in \cite{jig1}.  We now consider the case that $1< p<\infty$ for   type 1 subdiagonal algebras.

  Let $1< p,q<\infty$  such that $\frac1p +\frac1q=1$. It is known that there exists a bounded linear projection $F^p$ form  $L^p(\mathcal M)$  onto $H^p$ from \cite[Theorem 3.3]{jig1}. We also have that  the dual space of $H^p$ is
conjugate isomorphic to $H^q$(\cite[Corollary 3.4]{jig1}).  Thus for any $A\in \mathcal M$, we may define the right(resp. left) Toeplitz operator on $H^p$ as $t^p_Ax=F^p(xA)$(resp. $T^p_Ax=F^p(Ax)$) for all $x\in H^p$. It is elementary that $(t^p_A)^*=t^q_{A^*}$(resp. $(T^p_A)^*=T^q_{A^*}$) on $H^q$ for all $A\in\mathcal M$. If there no confusion to cause, then  we may simply  denote by $t_A$(resp. $T_A$) the right(resp. left)  Toeplitz operator defined by $A$ on $H^p$. If $A\in\mathfrak A$, then $t_A$(resp. $T_A$) is  said to be a right(resp.  left) analytic Toeplitz operator. We call the algebra $\mathcal R^p=\{t_A:A\in\mathfrak A\}$ and $\mathcal L^p=\{T_A:A\in\mathfrak A\}$ the right and left analytic Toeplitz algebras associated with $\mathfrak A$. It is interesting whether $\mathcal R^p$ and $\mathcal L^p$ are commutants of each other. For type 1 subdiagonal algebras, we have the following result. In fact,
we may consider more general case. We call a maximal subdiagonal algebra $\mathfrak A$ in $\mathcal M$ has  the weak factorization property if for any invertible operator $S\in \mathcal M$, there are a unitary operator $U\in \mathcal M$ and an operator $A\in\mathfrak A$ such that $S=UA$.  If $\mathfrak A$ has the universal factorization property, then  $\mathfrak A$ has  the weak one.

         {\bf Theorem 4.1.} Let   $\mathfrak A $  be a maximal    subdiagonal algebra with the weak factorization property.
 Then $\mathcal R^p$ and $\mathcal L^p$ are commutants of each other.

 \begin{proof}
  Note that $t_Ax=xA$ and $T_Ax=Ax$, $\forall x\in H^p$.
  It is known that $\mathcal L^{\prime}\supseteq  \mathcal R$.  Let $X\in\mathcal B(H^p)$ such that $XT_A=T_AX$ for all $A\in \mathfrak A$. Put $h=X(h_0^{\frac1p})\in H^p$. Then $X(Ah_0^{\frac1p})=XT_Ah_0^{\frac1p}=T_AX(h_0^{\frac1p})=Ah$, $\forall A\in\mathfrak A$. Take any $P\in\mathcal M_+$ and $\varepsilon>0$. Then $P+\varepsilon I\in\mathcal M_+$ is invertible and  there are a unitary $U\in\mathcal M$ and
   an operator $A\in\mathfrak A$ such that $P+\epsilon I=UA$ by  assumption. Thus
\begin{align*}
&\ \ \ \ \|(P+\epsilon I)h\|_p=\|UAh\|_p=\|Ah\|_p
=\|XAh_0^{\frac1p}\|_p\\
&\leq \|X\|\|Ah_0^{\frac1p}\|_p=\|X\|\|(P+\epsilon I)h_0^{\frac1p}\|_p .
\end{align*}
 It follows that $\|Ph\|_p\leq \|X\|\|Ph_0^{\frac1p}\|_p$ for all positive $P\in \mathcal M$.
 Now for any $B\in\mathcal M$, it is trivial that $\|B\xi\|_p=\||B|\xi\|_p$ for all $\xi\in L^p(\mathcal M)$.
   Then $\|Bh\|_p\leq \|X\|\|Bh^{\frac1p}_0\|_p$. We may define a bounded linear operator $Y$ on $  L^p(\mathcal M)$
   by $YBh_0^{\frac1p}=Bh$ for all $B\in\mathcal M$.  Note that for all $B,C\in\mathcal M$,
   $$L_CY(Bh_0^{\frac1p})=CBh=Y(CBh_0^{\frac1p}=YL_C(Bh_0^{\frac1p}).$$
   It follows that
    there is an element $G\in \mathcal M$ such that $Y=R_G$ from \cite[Corollary 1.6]{js}. Thus
$$XT_Ah_0^{\frac1p}=Ah=Y(Ah_0^{\frac1p})=Ah_0^{\frac1p}G$$ for all $A\in\mathfrak A$.
This implies that $G\in \mathfrak A$ by \cite[Theorem 2.7]{jig} and $X=t_G$. Therefore $(\mathcal L^p)^{\prime}=\mathcal R^p$. Similarly, we have $(\mathcal R^p)^{\prime}=\mathcal L^p$.  \end{proof}

   If   $\mathcal M$  is finite  or $\mathfrak A$ is of type 1, then   $\mathfrak A$ has the universal factorization property from \cite[Corollary 3.6]{jis} and Corollary 2.2. On the other hand, if we only consider analytic Toeplitz algebras, then Theorem 4.1 holds even if $p=1$.
To consider the reflexivity of analytic Toeplitz algebras, we summarize   the following elementary facts.

{\bf Proposition 4.2.} Let  $\mathfrak A$ be a type 1 subdiagonal algebra and $1\leq p<\infty$.

$(1)$    Suppose that  $x\in L^p(\mathcal M)$ is   a right separating and cyclic vector  such that $Ax=h_0^{\frac1p} $ for an operator  $A\in \mathcal M$. Then there is a contraction  $U\in\mathcal M$  and a $\xi\in H^p$ such that $x=U\xi $.

$(2)$  For any $\xi_1,\xi_2\cdots,\xi_n\in H^p$, there exists $\xi\in H^p$ such that $\|\xi_kA\|_p\leq \|\xi A\|_p $ for every $1\leq k\leq n$ and $A\in\mathfrak A$.

$(3)$ $H^p_0=\vee\{h_0^{\frac1p}U_{i_1}U_{i_2}\cdots U_{i_n}D:D\in\mathfrak D, U_{i_k}\in\mathcal U, n\geq 1\}$ and $H^p=L^P(\mathfrak D)+H_0^p$.

 \begin{proof}
 $(1)$  If $x=v|x|$ is the polar decomposition of $x$, then $V$ is unitary. We then may assume that $x\geq 0$. If $p=2$, then we easily have the result from Lemma 3.4 since $(h_0^{\frac12}A)x=h_0$. In fact, we have $U$ is unitary and $\xi $ is right outer.

 $(i)$  Let $2<p<\infty$. Then $\frac1p +\frac1r=\frac12 $ for some  $r>1$. Note that  $y=xh^{\frac1r}_0\in L^2(\mathcal M)$ is a separating and cyclic vector such that $Ay=h_0^{\frac12}$. Then  $y=U\eta$ for  a unitary operator $U$ and a right  outer element $\eta\in H^2$ as above.
 This means that $xh_0^{\frac1r}=U\eta$ and thus $\eta=U^*xh_0^{\frac1r}$.  We claim that $\xi=U^*x\in H^p$. In fact,
 $tr(\eta h_0^{\frac12} B)=tr(U^*xh_0^{\frac1r}h_0^{\frac12}B)=tr(\eta h_0^{\frac1q}B) =0$ for all $B\in\mathfrak A_0$.  Then $\xi\in H^p$ and $x=U\xi$.

$(ii)$ $1\leq p<2$. Put  $x=x^{\frac{p}{2}}x^{\frac{p}{r}}$, where $\frac1p=\frac12+\frac1r$. Note that
 $x^{\frac{p}r}\in   L^r(\mathcal M)$ and  $(Ax^{\frac{p}{2}})(x^{\frac{p}{r}}h_0^{\frac1q})=Axh_0^{\frac1q}=h_0$.
By Lemma 3.4, we have that $[x^{\frac{p}{r}}h_0^{\frac1q}\mathfrak A]_2=WH^2$ for a unitary operator $W$. Thus $x^{\frac{p}{r}}h_0^{\frac1q}=W\eta$ for some right outer element  $\eta\in H^2$. Then $\eta= W^*x^{\frac{p}{r}}h_0^{\frac1q}$.
This implies that $\eta_1=W^*x^{\frac{p}{r}}\in H^r$. Thus $x^{\frac{p}{r}}=W\eta_1$.

Now  $x^{\frac{p}2}W\in  L^2(\mathcal M)$ is a right separating and cyclic vector.
 We claim again that $[x^{\frac{p}2}W\mathfrak A]_2=VH^2$ for a partial isometry $V\in\mathcal M.$  Let $[x^{\frac{p}2}W\mathfrak A]_2=VH^2\oplus^{col} E  L^2(\mathcal M)$ such that $ E  L^2(\mathcal M)$ is of type 2 By Lemma 3.3. Then  $[AE L^2(\mathcal M) ]_2\subseteq [A ([x^{\frac{p}2}W\mathfrak A]_2) ]_2 $.  Note that $L^2(\mathcal M)H^2=L^1(\mathcal M)$.
    then $[AE L^1(\mathcal M) ]_1\subseteq [A ([x^{\frac{p}2}W\mathfrak A]_2)H^2]_1=[Axh_0^{\frac1q}\mathfrak A]_1=H^1$. It is known that
$[[AE L^1(\mathcal M) ]_1\mathfrak A_0]_1=[AE L^1(\mathcal M) ]_1$ since $EL^2(\mathcal M)$  is of type 2.  Similar to the proof in Lemma 3.4, we have that
$E=0$.   Note that   $V$ is a co-isometry since $ x^{\frac{p}2}W=Vy $ for some $y\in H^2$.    Thus $x^{\frac{p}{2}}W=V\xi_1$ for $\xi_1\in H^2$. Put  $ \xi=\xi_1\eta_1\in H^p$. $x=V\xi$.

$(2)$ Let $x=(\sum\limits_{k=1}^m|\xi|^2 +h_0^{\frac2p})^{\frac12}\in  L^p(\mathcal M).$ Then $\xi^*\xi \leq x^2$ as well as $(h_0^{\frac1p})^2\leq x^2$. By \cite[formula(1.1)]{js},  there are contractions  $B, B_1,\cdots,B_m\in\mathcal M$ such that
$h_0^{\frac1p}=Bx$ and $\xi_k=B_kx$ for $1\leq k\leq m$.   We  then have a contraction $C\in \mathcal M$ and $\xi\in H^p$  such that $x=B\xi$ as above.  Thus $\|\xi_k A\|_p=\|B_k x A\|_p=\|B_kC\xi A\|_p\leq\|\xi A\|_p$  for every $1\leq k\leq m$ and $A\in\mathfrak A$.

$(3)$  This  is  elementary  by Theorem 2.5. \end{proof}

 {\bf Theorem 4.3.} Let $\mathfrak A$ be a type  1  subdiagonal algebra   with multiplicity 1   and $1<p<\infty$.  Then the right analytic Toeplitz operator algebra $\mathcal R^p$ on $H^p$ is hereditary  reflexive.

 \begin{proof}
   By \cite[Proposition 4.3]{pel} and  Proposition 4.2(2), it is sufficient to show that $\mathcal R^p$ is reflexive.
 Let $U\in\mathcal M$ be a partial isometry with $U^*U\in\mathfrak D$  such that $H_0^2=UH^2$ .
 Then $\mathfrak A= \bigvee\{U^m\mathfrak D: m\geq 0\}$ by Theorem 2.5.

   Let $A\in algLat\mathcal R^p$. Then $A^*\in algLat\mathcal (R^p)^*$.  Moreover,  $EH^q\in Lat \mathcal R_p^*$  for any projection $E\in \mathfrak D$. Thus $A^*EH^q\subseteq EH^q$ for all projection $E\in\mathfrak D$. It easily follows that $A^*T_E=T_EA^*$ and thus $A^*T_D=T_DA^*$ for all $D\in\mathfrak D$.

   For any $z\in\mathbb D=\{z\in\mathbb C:|z|<1\}$, put $K_z=\{(I-zU)^{-1}x: x\in L^q(\mathfrak D)\}$.
 We claim that $K_z\in Lat\mathcal (R^p)^*$  for all $z\in\mathbb D$. In fact,  it is trivial that $K_z$ is closed and  $K_z \mathfrak D\subseteq K_z$.  It is sufficient  to show that $t_{U^*}K_z\subseteq K_z$.

 Let $\mathcal T$ be the set of entire  elements of
 $\mathcal M$,  that is,  those elements  $X\in\mathcal M$ for which the function
 $t\to\sigma_t^{\varphi}(X)$ can be extended to  an $\mathcal M$-valued entire function over
 $\mathbb C$. For any $X\in\mathcal M$ and $r>0$, we let
 \begin{equation}X_r=\sqrt{\frac{r}{\pi}}\int_{\mathbb
 R}e^{-rt^2}\sigma_t^{\varphi}(X)dt. \end{equation} Then
 $X_r\in\mathcal T$ and
 \begin{equation} \sigma_z^{\varphi}(X_r)=\sqrt{\frac{r}{\pi}}\int_{\mathbb
 R}e^{-r(t-z)^2}\sigma_t^{\varphi}(X)dt ,\ \ \forall z\in\mathbb C \end{equation} by the proof of Lemma VIII 2.3
 in  \cite{tak1}.
  Note that  $\sigma_{z}^{\varphi}(X)=h_0^{iz}Xh_0^{-iz}$ for any $z\in \mathbb
 C$ and  any entire element $X$ from $(2.1)$.  By Lemma VI 2.4 in    \cite{tak1},
 $X_r$ converges $\sigma$-weakly to
$X$ as $r\to \infty$.
Now replacing $X$ by $U$, we have
\begin{equation}U_rh_0^{\frac1q}=h_0^{\frac1q}\sigma_{\frac{i}{q}}(U_r).\end{equation} However, By formula $(4.2)$,
\begin{align*} \sigma_{\frac{i}{q}}^{\varphi}(U_r)&=\sqrt{\frac{r}{\pi}}\int_{\mathbb
 R}e^{-r(t-\frac{i}q)^2}\sigma_t^{\varphi}(U)dt \\
 &=\sqrt{\frac{r}{\pi}}\int_{\mathbb
 R}e^{-r(t-\frac{i}q)^2}UD(t)dt\\
 &=U\sqrt{\frac{r}{\pi}}\int_{\mathbb
 R}e^{-r(t-\frac{i}q)^2}D(t)dt, \end{align*}
 where $D (t)=U^*\sigma_t^{\varphi}(U )\in\mathfrak D$  by formula $(2.6)$.
  Note that
    \begin{equation}  \sqrt{\frac{r}{\pi}}\int_{\mathbb
 R}e^{-r(t-\frac{i}q)^2}D(t)dt=\sqrt{\frac{r}{\pi}}\int_{\mathbb
 R}e^{-r(t-\frac{i}q)^2}U^*\sigma_t^{\varphi}(U)dt
  =U^*\sigma_{\frac{i}{q}}^{\varphi}(U_r)\in \mathfrak D.
 \end{equation}
  It is known that $\sigma_{\frac{i}{q}}^{\varphi}(U_r)$ is bounded.  Without loss of generality, we may assume that $\sigma_{\frac{i}{q}}^{\varphi}(U_r)$ $\sigma$-weakly to   $V\in\mathfrak A_0$ as $r\to \infty.$ We now have $U^*V\in\mathfrak D$ as above $(4.4)$  and  $V=UU^*V$. It follows from $(4.3)$ that
 $Uh_0^{\frac1q}= h_0^{\frac1q}UU^*V$.  Then $h_0^{\frac1q}U^*=(V^*U)U^*h_0^{\frac12}$.

 It trivial that $t_{U^*}D h_0^{\frac1q}=0$ for all $D\in\mathfrak D$.
 \begin{align*}
 \ \ \ \ t_{U^*} (U^mDh_0^{\frac1q})
 &=F^q(U^m Dh_0^{\frac1q}U^*)
 =F^q( U^{m-1}U D(h_0^{\frac1q}U*))\\
 &= U^{m-1}  UD(V^*U)U^*h_0^{\frac1q}
 =U^{m-1}D_1 h_0^{\frac1q}
 \end{align*} for any $m\geq 1$,
  where $D_1= U(DV^*U)U^*\in\mathfrak D$ by Proposition 2.3.  Thus $t_{U^*}K_z\subseteq K_z$ for all $z\in\mathbb D$. Hence $K_z\in Lat(\mathcal R^p)^*$.
 Then $A^*K_z\subseteq K_z$ for all $z\in\mathbb D$.

 On the other hand, for any $\xi=(I-zU)^{-1}x=\sum\limits_{n=0}^{\infty}z^n U^nx\in K_z$,  $A^*\xi\in K_z$. $T_{U^*}\xi=\sum\limits_{n=1}^{\infty}z^n U^*U^nx=zU^*U\xi$. Thus $T_{U^*}A^*\xi=zU^*UA^*\xi$. Moreover, $A^*T_{U^*}\xi=zA^*U^*U\xi=zA^*T_{U^*U}\xi=T_{U^*U}A^*\xi=U^*UA^*\xi$ since $U^*U\in\mathfrak D$. That is,
 $A^*T_{U^*}\xi=T_{U^*}A^*\xi$ for all $\xi\in K_z$ and $z\in\mathbb D$.

 We next claim that $\bigvee\{K_z:z\in\mathbb D\}=H^q$. In fact, the dual space of $H^q$ is conjugate isomorphic to $H^p$ by \cite[Corollary 3.4]{jig1}. Take any  $\eta\in H^p$ such that
 $$\langle(I-zU)^{-1}x,\eta\rangle =\langle \sum\limits_{m\geq 0}z^mU^mx, \eta\rangle=\sum\limits_{m\geq 0}z^m\langle U^mx,\eta\rangle=0$$
 for all $x\in L^q(\mathfrak D)$ and $z\in \mathbb D$. Then $\langle U^mx,\eta\rangle=0 $ for all $x\in L^q(\mathfrak D)$ and $m\geq 0$. It follows that $\eta=0$ by Proposition 4.2(3).
 Then $A^*T_{U^*}=T_{U^*}A^*$. Hence
 $AT_U=T_UA$. Thus $A\in(\mathcal L^p)^{\prime}=\mathcal R^p $ by Theorem 4.1. That is, $\mathcal R^p$ is reflexive.
 Consequently, $\mathcal R^p$ is hereditary  reflexive by \cite[Proposition 4.3]{pel} and  Proposition 4.2(2).  \end{proof}

{\bf Corollary 4.4.} Let  $\mathfrak A=\mathfrak D\rtimes_{\alpha} \mathbb Z_+ $  be a non-self-adjoint  crossed product.  Then  the  right analytic Toeplitz algebra $\mathcal R^p$ on $H^p$($1<p<\infty$)  associated with $\mathfrak A$  is hereditary reflexive.

%\section*{References}

\end{document}